\documentclass[10pt]{article}
\usepackage[top=1in, bottom=1in, left=0.75in, right=0.75in]{geometry}

\usepackage{caption}
\usepackage{cite}
\usepackage{amsmath,amssymb,amsfonts}
\usepackage{algorithmic}
\usepackage{graphicx}
\usepackage{textcomp}
\def\BibTeX{{\rm B\kern-.05em{\sc i\kern-.025em b}\kern-.08em
    T\kern-.1667em\lower.7ex\hbox{E}\kern-.125emX}}

\usepackage[utf8]{inputenc}
\usepackage{listings}
\usepackage{xspace}
\usepackage{soul}
\usepackage{float}
\usepackage{booktabs}
\usepackage{multicol}
\usepackage{wrapfig}

\newcommand{\by}{\ensuremath{\mathbf{y}}\xspace}

\newcommand{\bx}{\ensuremath{\mathbf{x}}\xspace}

\newcommand{\bxdot}{\ensuremath{\mathbf{\dot{x}}}\xspace}

\newcommand{\bxhat}{\ensuremath{\mathbf{\hat{x}}}\xspace}

\newcommand{\bxdothat}{\ensuremath{\mathbf{\hat{\dot{x}}}}\xspace}

\title{\bf Numerical differentiation of noisy data:\\ A unifying multi-objective optimization framework}
\author{Floris van Breugel$^1$,
J. Nathan Kutz$^2$,
Bingni W. Brunton$^3$}
\date{}

\begin{document}

\maketitle

\footnotesize{$^1$Department of Mechanical Engineering, University of Nevada, Reno, NV 89557 (e-mail: fvanbreugel@unr.edu)}

\footnotesize{$^2$Department of Applied Math, University of Washington, Seattle, WA, 98195, USA (e-mail: kutz@uw.edu)}

\footnotesize{$^3$Department of Biology, University of Washington, Seattle, WA, 98195 (e-mail: bbrunton@uw.edu)}



\begin{abstract}
Computing derivatives of noisy measurement data is ubiquitous in the physical, engineering, and biological sciences, and it is often a critical step in developing dynamic models or designing control. 
Unfortunately, the mathematical formulation of numerical differentiation is typically ill-posed, and researchers often resort to an \textit{ad hoc} process for choosing one of many computational methods and its parameters.
In this work, we take a principled approach and propose a multi-objective optimization framework for choosing parameters that minimize a loss function to balance the faithfulness and smoothness of the derivative estimate. 
Our framework has three significant advantages. First, the task of selecting multiple parameters is reduced to choosing a single hyper-parameter.  
Second, where ground-truth data is unknown, we provide a heuristic for automatically selecting this hyper-parameter based on the power spectrum and temporal resolution of the data.  
Third, the optimal value of the hyper-parameter is consistent across different differentiation methods, thus our approach unifies vastly different numerical differentiation methods and facilitates unbiased comparison of their results. 
Finally, we provide an extensive open-source Python library \texttt{pynumdiff} to facilitate easy application to diverse datasets (\texttt{https://github.com/florisvb/PyNumDiff}). 

\end{abstract}

\vspace{10pt}
\emph{Keywords:} Numerical differentiation, derivatives, optimization, data-driven modeling.

\normalsize{}
\vspace{10pt}
\begin{multicols}{2}
\section{Introduction} \label{sec:intro}
Derivatives describe many meaningful characteristics of physical and biological systems, including spatial gradients and time rates-of change.
However, these critical quantities are often not directly measurable by sensors.
Although computing derivatives of analytic equations is straightforward, estimating derivatives from real sensor data remains a significant challenge because sensor data is invariably corrupted by noise~\cite{Ahnert2007}.
More accurate estimation of derivatives would improve our ability to produce robust diagnostics, formulate accurate forecasts, build dynamic or statistical models, implement control protocols, and inform policy making.
There exists a large and diverse set of mathematical tools for estimating derivatives of noisy data, most of which are formulated as an ill-posed problem regularized by some appropriate smoothing constraints.  
However, the level and type of regularization are typically imposed in an {\em ad hoc} fashion, so that there is currently no consensus ``best-method" for producing ``best-fit" derivatives.

One particularly impactful application of estimating derivatives is the use of time-series data in modeling complex dynamical systems.
These models are of the form $d{\bx}/dt = \dot{{\bx}} = f({\bx})$, where ${\bx}$ is the state of the system.
Models of this kind have been integral to much of our understanding across science and engineering~\cite{lin1988mathematics}, including in classical mechanics~\cite{goldstein2002classical}, electromagnetism~\cite{jackson2007classical}, quantum mechanics~\cite{griffiths2018introduction}, chemical kinetics~\cite{masel2001chemical}, ecology~\cite{kot2001elements} epidemiology~\cite{rothman2008modern}, and neuroscience~\cite{Roth2014, Madhav2020, keel2010generalized}.
In some cases, even higher order time derivatives are also crucial for understanding the dynamics~\cite{Lin2019}.
A recent innovation in understanding complex dynamical systems uses data-driven modeling, where the underlying dynamics are learned directly from sensor data using a variety of modern methods~\cite{Schmidt2009,Brunton2016pnas,Daniels2015}.
For this application in particular, a derivative with both small and unbiased errors is crucial for learning interpretable dynamics.

\begin{figure*}[t] 
\centering
\includegraphics[width=2\columnwidth]{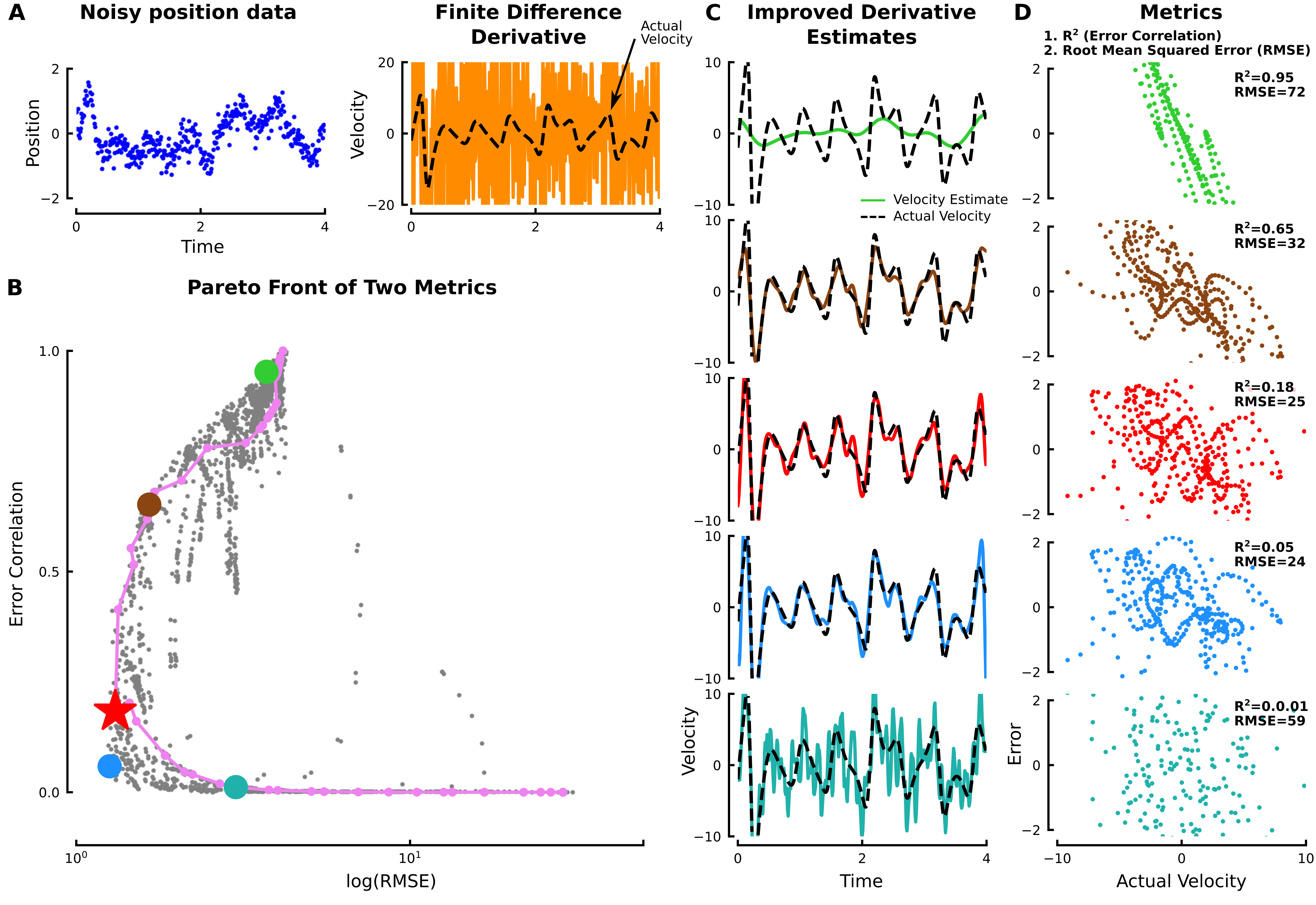}
\caption{
Choice of parameters leads to a diversity of derivative estimates. A. Noisy time series data, from a Lorenz system, and the corresponding finite difference derivative. 
B. To evaluate the quality of a derivative estimate relative to the ground truth, we consider two metrics: Root Mean Square Error (RMSE), and the Pearson's correlation coefficient ($R^2$) between the error and the true value of the derivative. Gray dots show the values of these metrics for $5,481$ different sets of parameter choices for a smoothed Savitzky-Golay filter. The violet line is the result of our multi-objective optimization framework and nearly traces the Pareto front of the metrics. The derivative estimates and metrics for the five colored points along the Pareto front are shown in C and D, respectively.}
\label{fig:fig_1}
\end{figure*}

In principle, the discrete derivative of position can be estimated as the finite difference between adjacent measurements.
If we write the vector of all noiseless positions in time measured with timestep $\Delta t$ as $\bx$, then
\begin{equation} \label{eq:finite_difference}
  \bxdot_k = \frac{\bx_{k+1}-\bx_k}{\Delta t},
\end{equation}
where $k$ indexes snapshots in time.
In reality, however, only noisy measurements $\by$ are available, 
\begin{equation*}
    \by=\bx+\boldsymbol\eta,
\end{equation*}
where $\boldsymbol\eta$ represents measurement noise.
Here we will assume $\boldsymbol\eta$ is zero-mean Gaussian noise with unknown variance.
Even with noise of moderate amplitude, a na\"{i}ve application of Eq.~\eqref{eq:finite_difference} produces derivative estimates that are far too noisy to be useful (Fig.~\ref{fig:fig_1}A). 
Thus, more sophisticated methods for data smoothing and/or differentiation of noisy time series measurements of position $\by$ are required. 

Although smoothing mitigates  the errors, it can also introduce biases. 
Our goal in this paper is to develop a general approach for methodically choosing parameters that balance the need to minimize both error and bias. 
We use $\bxhat$ and $\bxdothat$ to denote the smoothed \textit{estimates} of the position and its derivative computed from $\by$, respectively. 
To evaluate the quality of these estimates, we compare these estimates to the \textit{true} discrete time position and its derivative, $\bx$ and $\bxdot$. 
Developing approaches for estimating $\bxdothat$ from noisy measurements $\by$ has been the focus of intense research for many decades. 
Despite the diversity of methods that have been developed, only a few studies have performed a comprehensive comparison of their performance on different types of problems~\cite{Walker1998, Crenshaw2000, Ahnert2007}.

In this paper, we tackle the challenge of parameter selection by developing a novel, multi-objective optimization framework for choosing parameters to estimate the derivative of noisy data. 
Our approach minimizes a loss function consisting of a weighted sum of two metrics computed from the derivative estimate: the faithfulness of the integral of the derivative and its smoothness.
We suggest these metrics as proxies for minimizing the error and bias of the estimated derivative, and we show that sweeping through values of a single hyper-parameter $\gamma$ produces derivative estimates that generally trace the Pareto front of solutions that minimize error and bias.
Importantly, this optimization framework assumes no knowledge of the underlying true derivative and reduces the task of selecting many parameters of any differentiation algorithm to  solving a loss function with a single hyper-parameter. 
Furthermore, we show that the value of the hyper-parameter is nearly universal across four different differentiation methods, making it possible to compare the results in a fair and unbiased way. 
For real-world applications, we provide a simple heuristic for automatically determining a value of $\gamma$ that is derived from the power spectrum and temporal resolution of the data.
All of the functionality described in this paper is implemented in an open-source Python toolkit \texttt{pynumdiff}, which is found here: \texttt{https://github.com/florisvb/PyNumDiff}.

\section{Motivation for error metrics} \label{sec:motivation}

What is a ``good'' estimate of a derivative? 
Let us start by considering a toy system with synthetic measurement noise, where we are able to evaluate the quality of an estimated derivative by comparing to the true, known derivative.
We consider two metrics for evaluating the quality of a derivative Fig.~\ref{fig:fig_1}B--D; later, we use these same metrics to evaluate the performance of our optimization framework, which does not have access to the ground truth. 

First, the most intuitive metric is how faithfully the estimated derivative $\bxdothat$ approximates the actual derivative $\bxdot$. 
We can measure this using the root-mean-squared error, 
\begin{equation} \label{eq:RMSE}
  \text{RMSE}(\bxdothat, \bxdot) = \lVert(\bxdothat-\bxdot)\rVert_2,
\end{equation}
where $\lVert \cdot \rVert_2$ is the vector 2-norm. 
If the data are very noisy, a small RMSE can only be achieved by applying significant smoothing. 
However, smoothing the data often attenuates sharp peaks in the data and results in underestimating the magnitude of the derivative.

To measure the degree to which the derivative estimate is biased due to underestimates of the actual derivative, we calculate the square of the Pearson's correlation coefficient, $R^2$, between the errors $(\bxdothat-\bxdot)$ and the actual derivative $\bxdot$. 
We refer to this metric as the \textit{error correlation}, which is bounded between 0 and 1.
Small error correlations imply that the imposed dynamics of the differentiation  method (e.g. filtering) minimally influenced the derivative estimate; therefore, the method of estimating derivatives would have minimal impact on any models that are constructed using these estimates. 
Conversely, large error correlations imply that the estimate is significantly influenced by the dynamics of the differentiation method and typically correspond to very smooth estimates. 
In the limit where the derivative estimate is a horizontal line, the error correlation takes on a value of unity.
Other metrics that measure the smoothness, for example the total variation or tortuosity, may be substituted for error correlation~\cite{Ahnert2007}; however, these metrics are harder to interpret.
For instance, if the true derivative is very smooth, a low total variation is desired, whereas if the true derivative is quite variable, a high total variation would correspond to an accurate derivative.  
In contrast, a low error correlation is desirable for any true derivative.  

For many datasets, the RMSE and error correlation metrics define a Pareto front, where no single parameter choice minimizes both values (Fig.~\ref{fig:fig_1}B). 
Furthermore, the minimal RMSE can be achieved with a variety of different error correlations. 
The most suitable parameter set depends on the application of the estimated derivative: is a non-smooth derivative with minimal bias preferred (Fig.~\ref{fig:fig_1}B-D: teal), or one that is smooth, but biased (Fig.~\ref{fig:fig_1}B-D: brown).
We suggest that, for most purposes, the estimated derivative that balances these metrics (Fig.~\ref{fig:fig_1}B-D: blue and red) serves as a reasonable starting point. 

\section{Methods for Numerical Differentiation} \label{sec:diff_methods}

A large variety of methods for numerical differentiation exist, and a complete review of them all is beyond the scope of this paper. 
Instead, we have selected four differentiation methods (Table~\ref{table:diffmethods}), which make different assumptions and represent different approaches to computing the derivative including both global and local methods \cite{Ahnert2007}, to showcase the universal application of our optimization framework. 

One common approach to manage noisy data is to apply a smoothing filter to the data itself, followed by a finite difference calculation. 
In this family of differentiation methods, we chose to highlight the \textbf{Butterworth filter} \cite{Butterworth1930}, which is a global spectral method with two parameters: filter order and frequency cutoff. 

The second family of methods relies instead on building a local model of the data through linear regression. 
A common and effective approach involves making a sliding polynomial fit of the data~\cite{Belytschko1996}, often referred to as locally estimated scatterplot smoothing (LOESS)~\cite{harrell2015regression}. 
An efficient approach for accomplishing the same calculations is the \textbf{Savitzky-Golay filter}, which builds the polynomial model in the frequency domain~\cite{Schafer2011, SavGol}. 
The Savitzky-Golay filter has two parameters: window size and polynomial order. 
By default, a Savitzky-Golay filter provides a jagged derivative because the polynomial models can change from one window to the next, so here we also apply some smoothing by convolving the result with a Gaussian kernel.
This smoothing adds a third parameter: a smoothing window size. 

The third family we consider is the \textbf{Kalman filter}~\cite{Kalman1960, zarchan2015fundamentals, Aravkin2017}. 
The Kalman filter is most effective when models of the system and of the noise characteristics are known. 
Our focus here is the case where neither is known, so we chose to highlight a constant acceleration forward-backward Kalman smoother~\cite{crassidis_junkins} with two parameters: the model and noise covariances. 

Finally, we consider an optimization approach to computing derivatives with the \textbf{total variation regularization} (TVR) method~\cite{OsherFatemi92, Chartrand11}. 
One advantage of the TVR methods is that there is only a single parameter, which corresponds to the smoothness of the derivative estimate. 
TVR derivatives are not as widely used as the other three methods we highlight, so we provide a brief overview here. 
Solving for the TVR derivative involves first finding $\bxhat$ and its corresponding finite-difference derivative $\bxdothat$ (calculated according to Eq.~\ref{eq:finite_difference}) that minimize the following loss function, 
\begin{equation}
L = \lVert \by - \bxhat \rVert_2 + \gamma*TV(\bxdothat).
\label{eq:totalvariationregularization}
\end{equation}
Here $TV$ is the total variation,
\begin{equation}
TV(\bxdothat) = \frac{1}{m}\left\lVert\bxdothat_{0:m-1}-\bxdothat_{1:m}\right\rVert_{1},
\label{eq:totalvariation}
\end{equation}
where $\lVert\cdot\rVert_1$ denotes the $\ell_1$ norm and $m$ is the number of time snapshots in the data. 
The single parameter for this method is $\gamma$, and larger values result in smoother derivatives. 
If $\gamma$ is zero, this formulation reduces to a finite difference derivative. 

Solutions for TVR $\bxdothat$ can be found with an iterative solver~\cite{Chartrand11}. 
Because both components of the loss function Eq.~\eqref{eq:totalvariationregularization} are convex, we can also solve for $\bxdothat$ using convex optimization tools, such as cvxpy \cite{cvxpy}, and with a convex solver, such as MOSEK \cite{mosek}. 
The two methods are equivalent, if the iterative solver is repeated sufficiently many times. 

The convex solution to penalizing the first order difference in time, as in Eq.~\eqref{eq:totalvariation}, results in a piece-wise constant derivative estimate. 
By offloading the calculations to a convex optimization solver, however, we can easily penalize higher order derivatives by replacing the 1$^{st}$ order finite difference derivative $\bxdothat$ in Eq.~ \eqref{eq:totalvariationregularization} with a 2$^{nd}$ order ($\mathbf{\hat{\ddot{x}}}$) or 3$^{rd}$ order ($\mathbf{\hat{\dddot{x}}}$) finite difference derivative. 
Penalizing higher-order time derivatives results in smoother derivative estimates. 
For example, penalizing the 2$^{nd}$ order derivative results in a piece-wise linear derivative estimate, whereas penalizing the 3$^{rd}$ order derivative, also known as the \textit{jerk}, results in a smooth estimate. 
In this paper, we will use the total variation regularized on the jerk (TVRJ). 
For large datasets, solving for the TVRJ derivative is both computationally expensive and can accumulate small errors. 
To manage the size of the optimization problem, we solve for the TVRJ derivative in sliding windows.

\begin{table*}[t]
\footnotesize
\centering
\begin{tabular}{@{}lcccc@{}}
\toprule
\multicolumn{1}{c}{\textbf{Full name}}                 & \textbf{\begin{tabular}[c]{@{}c@{}}Abbreviated\\ name\end{tabular}} & \textbf{\# Parameters} & \textbf{\begin{tabular}[c]{@{}c@{}}Computational\\ cost\end{tabular}} & \textbf{References} \\ \midrule
Butterworth filter followed by finite difference       & Butterworth                                                         & 2                      & low                                                                   & \cite{Butterworth1930}                 \\
Smooth Savitzky-Golay filter                           & Savitzky-Golay                                                      & 3                      & low                                                                   & \cite{Schafer2011, SavGol}                 \\
Constant acceleration forward-backward Kalman smoother & Kalman smooth                                                       & 2                      & high                                                                  & \cite{crassidis_junkins}                 \\
Total Variation Regularized Jerk                       & TVRJ                                                                & 1                      & high                                                                  & \cite{Chartrand11}                 \\ \bottomrule
\end{tabular}
\caption{Summary of the four differentiation methods highlighted in this paper. \label{table:diffmethods}}
\end{table*}

\section{Computing derivatives of noisy data with no ground truth}

With noisy data collected in the real world, no ground truth is accessible.
The RMSE and error correlation metrics described in the previous section cannot be calculated and used to optimize parameter choices, so the parameter selection is an ill-posed problem. 
Even so---somehow---parameters must be chosen. 
In this section, we propose a general approach for choosing parameters and show that for a wide range of problems, noise levels, time resolutions, and methods, our approach yields reasonable derivative estimates without the need for hyper-parameter turning. 

\subsection{Optimization framework without ground truth derivatives}

Given noisy position measurements $\by$, we seek to estimate the derivative in time of the dynamical system that underlies the measurements $\bxdothat$.
When the ground truth $\bxdot$ is unknown, we propose choosing the set of parameters $\Phi$ (for any given numerical algorithm, including those enumerated in Table~\ref{table:diffmethods}) that minimize the following loss function, which is inspired by Eq.~\eqref{eq:totalvariationregularization},

\begin{equation}
L = \mbox{RMSE} \bigg( \mbox{trapz}(\bxdothat(\Phi)) + \mu, \by \bigg) + \gamma \bigg({TV}\big(\bxdothat(\Phi)\big)\bigg),
\label{eq:loss}
\end{equation}
where $\mbox{trapz}(\cdot)$ is the discrete-time trapezoidal numerical integral, $\mu$ resolves the unknown integration constant,
\begin{equation}
\mu=\frac{1}{m}\sum_{k=0}^{m}\bigg(\mbox{trapz}(\bxdothat(\Phi)) - \by\bigg),
\end{equation}
and $\gamma$ is a hyper-parameter. 
Note that this formulation has a single hyper-parameter $\gamma$, and a heuristic for choosing $\gamma$ is introduced in the following section. 

The first term of the loss function in Eq.~\eqref{eq:loss} promotes faithfulness of the derivative estimate by ensuring that the integral of the derivative estimate remains similar to the data, whereas the second term encourages smoothness of the derivative estimate. 
If $\gamma$ is zero, the loss function simply returns the finite difference derivative. 
Larger values of $\gamma$ will result in a smoother derivative estimate. 

This loss function effectively reduces the set of parameters $\Phi$ (which ranges between 1 and 3 or more, depending on the method) to a single hyper-parameter $\gamma$. 
Unfortunately, $L$ is not convex, but tractable optimization routines can be used to solve for the set of $\Phi$ that minimize $L$.
Here we use the Nelder-Mead method~\cite{Nelder1965}, a gradient descent algorithm, as implemented in SciPy~\cite{2020SciPy-NMeth}, with multiple initial conditions.

\subsection{Heuristics for automated hyper-parameter tuning of $\gamma$}

The advantages of our loss function in Eq.~\eqref{eq:loss} are that it does not require any ground truth data, and it simplifies the process of choosing parameters by reducing all the parameters associated with any given method for differentiation to a single hyper-parameter $\gamma$ corresponding to the how smooth the resulting derivative should be. 
To understand the qualities of the derivative estimates resulting from parameters selected by our loss function, we begin by analyzing the derivative estimates of noisy sinusoidal curves using the Savitzky-Golay filter and return to our original metrics, RMSE and error correlation to evaluate the results. 

\begin{figure*}[t] 
\includegraphics[width=\textwidth]{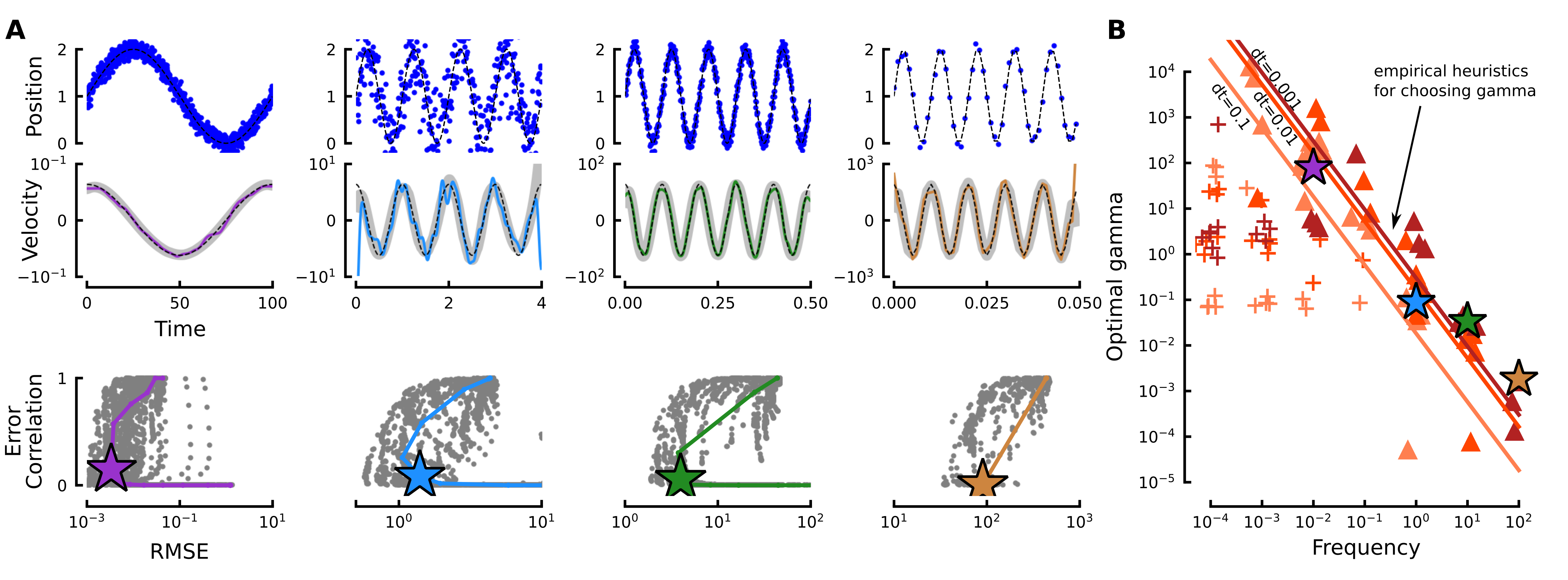}
\caption{Optimal choice of $\gamma$ is a function of frequency and temporal resolution of the data. A. (Top) Four example sine waves of different frequencies (note the time scales), temporal resolutions, and noise levels. (Middle) Comparison of the actual derivative (black dashed) with Savitzky-Golay estimates: lowest achievable RMSE (gray), and the result from our loss function with the optimal choice of $\gamma$ defined in the bottom panel. (Bottom) Trade-off between error correlation and RMSE for 5,481 potential parameter choices (gray) and the options provided by our loss function for a sweep through $\gamma$ (colored line). The star indicates the optimal choice of gamma, corresponding to the shoulder of the colored curve. (B) The optimal choice of $\gamma$ (defined in A) as a function of frequency (Hz), for different temporal resolutions of data (0.001, 0.01, 0.1 sec). Also included in the plot, but not indicated, are different noise levels (0-mean normally distributed with standard deviations of $0.05\%, 0.5\%, 5\%, $ and $25\%$ of the amplitude) and length of the dataset ($1, 4, 5, 25, 100, 500, 1000$ sec). The "+" markers indicate results from datasets for which the period was greater than the length of the time series, which were omitted from the fit. The diagonal lines indicate the empirical heuristics for choosing $\gamma$ based on a multivariate ordinary least squares model, provided in Eqn. \ref{eqn:goldgamma} and Table \ref{table:e2}.}
\label{fig:sinusoids}
\end{figure*}

Interestingly, sweeping through values of $\gamma$ results in derivative estimates with RMSE and error correlation values that generally follow the Pareto front defined by all possible derivative estimates for that given method (Fig.~\ref{fig:sinusoids}A). 
Which of these derivative estimates is best depends on the intended use of the derivative; nevertheless, we suggest that a good general purpose derivative is one that corresponds with the elbow in the lower left corner of the stereotypical curve traced by a sweep of $\gamma$ in the RMSE vs. error correlation space (the star-shaped markers in Fig.~\ref{fig:sinusoids}A). 
This point often, but not always, corresponds to the lowest RMSE (for example, see Fig.~\ref{fig:fig_1}). 
Although in many cases a quantitatively better derivative estimate than the one found by our loss function does exist (the gray dots in Fig.~\ref{fig:sinusoids}A that lie left of the star), the qualitative differences between these two derivative estimates are generally small (Fig.~\ref{fig:sinusoids}A middle row). 

{
\begin{wraptable}{c}{\columnwidth}
\footnotesize
\centering
\begin{tabular}{@{}lll@{}}
\toprule
\textbf{Variable} & \textbf{Coeff} & \textbf{P-value} \\ 
\midrule
intercept         & -5.26           & 0                \\
$\log(freq)$      & -1.55          & 0                \\
$\log(dt)$        & -0.74          & 0                \\
$\log(noise)$     & 0.11          & 0.32             \\
$\log(length)$    & 0.10           & 0.32            
\end{tabular}
\caption{Optimal $\log(\gamma)$ is correlated with frequency and temporal resolution, but not the noise or length of the dataset. The table provides the coefficients and associated p-values for a ordinary least squares model, with an adjusted $R^2=0.78$.
\label{table:e1}}
\vspace{10pt}

\begin{tabular}{@{}lll@{}}
\toprule
\textbf{Variable} & \textbf{Coeff} & \textbf{P-value} \\ \midrule
intercept         & -5.1           & 0                \\
$\log(freq)$      & -1.6          & 0                \\
$\log(dt)$        & -0.71          & 0               
\end{tabular}
\caption{Optimal $\log(\gamma)$ can be determined based on the frequency and temporal resolution of the data. The table provides the coefficients and associated p-values for a ordinary least squares model, with an adjusted $R^2=0.78$.\\
\label{table:e2}}
\vspace{5pt}
\end{wraptable}
}

In practice, the need to choose even a single parameter can be time consuming and arbitrary. 
To alleviate these issues, we derive an empirical heuristic to guide the choice of  $\gamma$ that corresponds with the elbow of the Pareto front. 
We found that the best choice of $\gamma$ is dependent on the frequency content of the data. 
To characterize this relationship, we evaluated the performance of derivative estimates achieved by a Savitzky-Golay filter by sweeping through different values of $\gamma$ for a suite of sinusoidal data with various frequencies ($f$), noise levels (additive white (zero-mean) Gaussian noise with variance $\sigma^2$), temporal resolutions ($\Delta t$), and dataset lengths (in time steps, $L$) (Fig.~\ref{fig:sinusoids}A-B). 

To describe this empirical relationship between the optimal choice of $\gamma$ and quantitative features of the data, we first considered an all-inclusive multivariate log-linear model,

\begin{equation}
\log(\gamma) = \alpha_1\log(f) + \alpha_2\log(\Delta t) + \alpha_3\log(\sigma) + \alpha_4\log(L) + \alpha_5.
\end{equation}

Fitting the data (Fig.~\ref{fig:sinusoids}B triangles) to this model with ordinary least squares resulted in an $R^2=0.76$, suggesting that, in many cases, it is feasible to automatically determine a reasonable guess for $\gamma$. 
Table~\ref{table:e1} provides the coefficients ($\alpha_k$) and associated p-values for each of the four terms and intercept. 
From this analysis we can conclude that the magnitude of measurement noise in the data is \textit{not} an important predictor of $\gamma$. 
We note, however, that here we have assumed that the magnitude of noise does not change within a time-series dataset.

Eliminating the unnecessary terms from our model results in slightly adjusted coefficients, provided in Table~\ref{table:e2}. 
In short, the optimal choice of $\gamma$, assuming that both low RMSE and low error correlation are valued, can be found according to the following relationship:
\begin{equation} \label{eqn:goldgamma}
\log(\gamma) = -1.6\log(f) -0.71\log(dt) -5.1.
\end{equation}

\begin{figure*}[p] 
\includegraphics[width=0.9\textwidth]{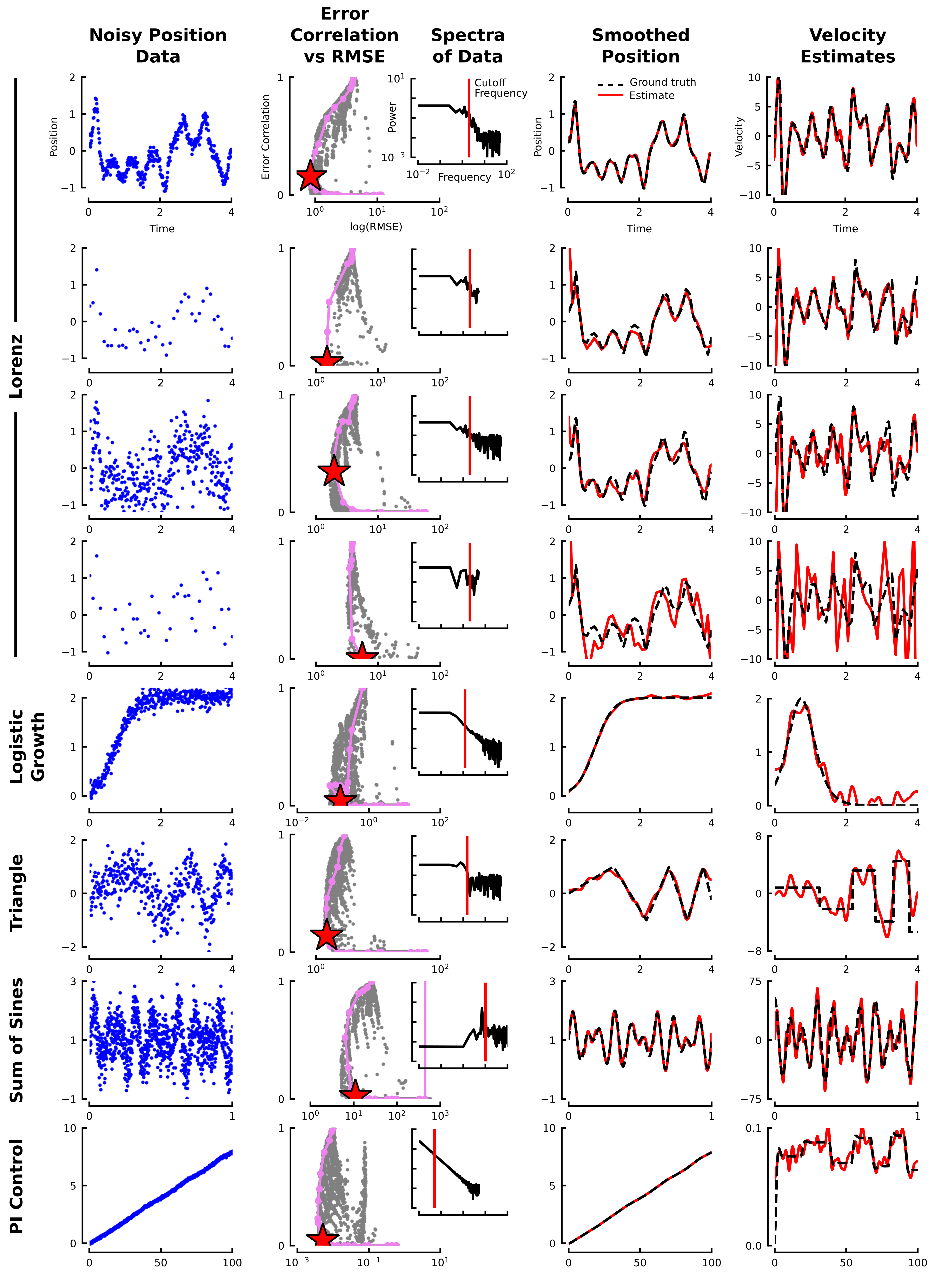}
\caption{Heuristic for choosing $\gamma$ is effective across a broad range of toy problems, using a Savitzky-Golay filter. The first column shows raw (synthetic) position data, indicating the shape of the data, degree of noise, and temporal resolution. Next we evaluate the performance of derivative estimate using the metrics described in the Fig. \ref{fig:fig_1}. Gray dots indicate the range of outcomes for 5,481 parameter choices, the violet line indicates the options provided by our loss function, and the red star indicates the performance using the automatically selected value of $\gamma$ according to Eqn. \ref{eqn:goldgamma}. Frequency of the data is evaluated by inspecting the power spectra; the red line indicates the frequency used to determine $\gamma$. The final two columns compare the ground truth and estimates for position and velocity.}
\label{fig:fig_3}
\end{figure*}

We analyze the performance of our loss function and heuristic with respect to a broad suite of representative synthetic problems. 
Real world data takes on a much greater diversity of shapes than the sinusoidal timeseries we used to derive the heuristic for choosing $\gamma$ given in Eq.~\eqref{eqn:goldgamma}. 
Because it is difficult to define a clear quantitative description of the range of shapes that real data might take on (such as frequency for a sinusoidal function), 
we first examine differentiating one component of a Lorenz system (Fig.~\ref{fig:fig_3}). 
From the power spectra, we select a frequency corresponding to the frequency where the power begins to decrease and the noise of the spectra increases. 
Although somewhat arbitrary, this approach (in conjunction with Eq.~\eqref{eqn:goldgamma}) allows us to use a standard signal processing tool to quickly determine a choice of $\gamma$. 
Our method produces reliable derivatives without further tuning in each case except high noise and low temporal resolution (Fig.~\ref{fig:fig_3}, fourth row), which is not surprising considering the low quality of the data.

Next we consider four other synthetic problems, all with similarly effective results (Fig.~\ref{fig:fig_3}). 
For the logistic growth problem, the curve traced by our loss function takes on a more complicated shape, perhaps because the characteristics of data vary substantially across time. 
Still, our heuristic results in a good choice of parameters that correspond to an accurate derivative. 
For the triangle wave, the loss function does a good job of tracing the Pareto front, and the heuristic selects an appropriate value of $\gamma$, yet the resulting derivative does show significant errors. This is likely due to two reasons. 
First, the Savitzky-Golay filter is designed to produce a smooth derivative, rather than a piece-wise constant one. 
Second, the frequency content of the data varies between two extremes, near-zero, and near-infinity. 
For the sum of sines problem, selecting the appropriate frequency cutoff is more straightforward than the previous problems, as we can simply choose a frequency shortly after the high frequency spike in the spectra. 
The final problem is a time-series resulting from a simulated dynamical system controlled by a proportional-integral controller subject to periodic disturbances. This data is a challenging problem for numerical differentiation, as the position data almost appears to be a straight line but does contain small variations. 
Our loss function does an excellent job of tracing the Pareto front in this case, and our heuristic results in an appropriate choice of $\gamma$. 

\begin{figure*}[htb] 
\includegraphics[width=0.9\textwidth]{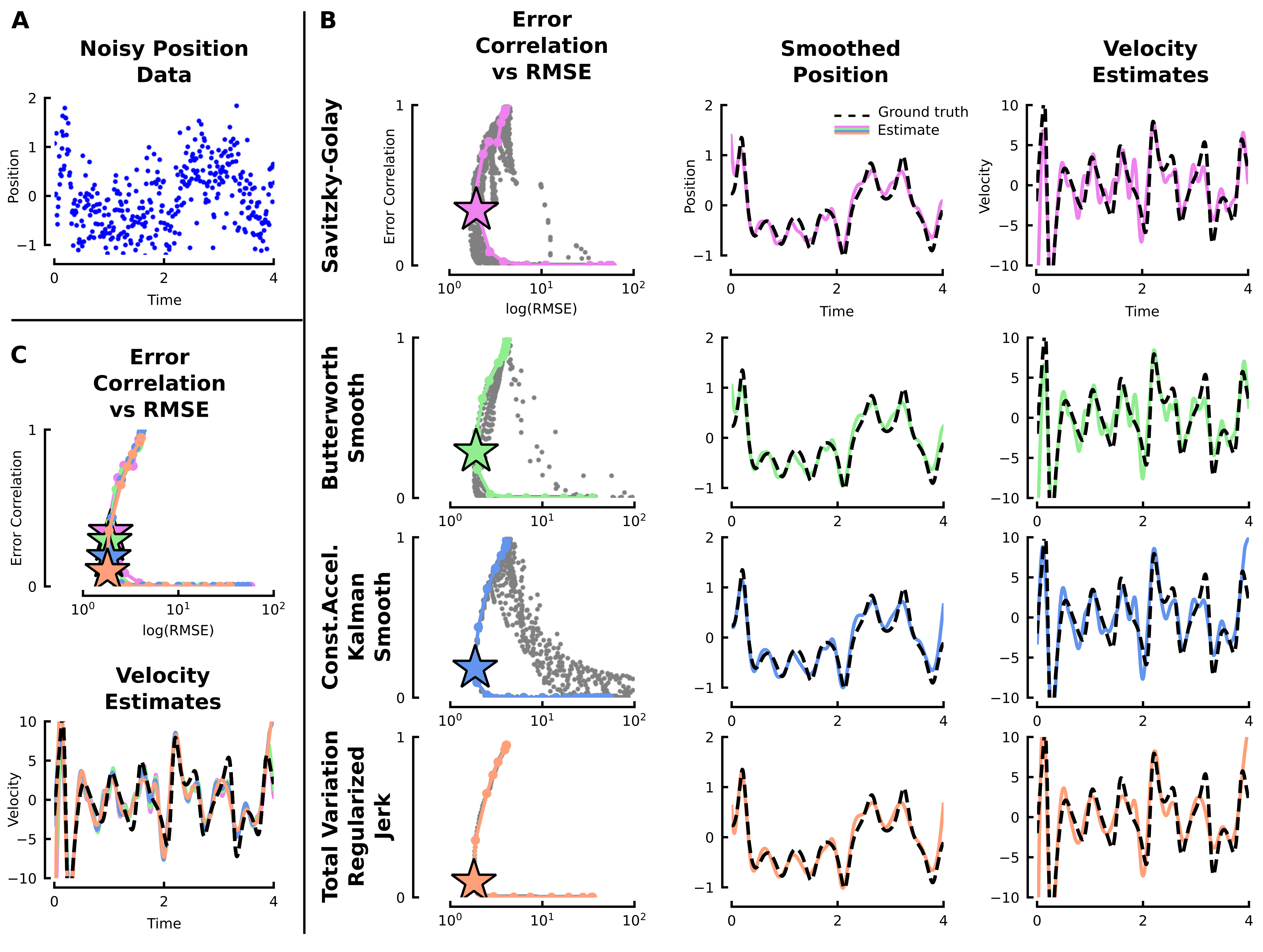}
\caption{Loss function and heuristic for choosing $\gamma$ is equally effective for different differentiation methods. A. Synthetic noisy data from the same Lorenz system as shown in Fig. \ref{fig:fig_3}. B. Comparison of metrics, position, and velocity estimates using four differentiation methods, with the same value of $\gamma$, as determined through the spectral analysis in Fig. \ref{fig:fig_3}. C. Overlay of the Pareto fronts and velocities for all four methods. }
\label{fig:fig_4}
\end{figure*}

\begin{figure*}[p] 
\includegraphics[width=0.9\textwidth]{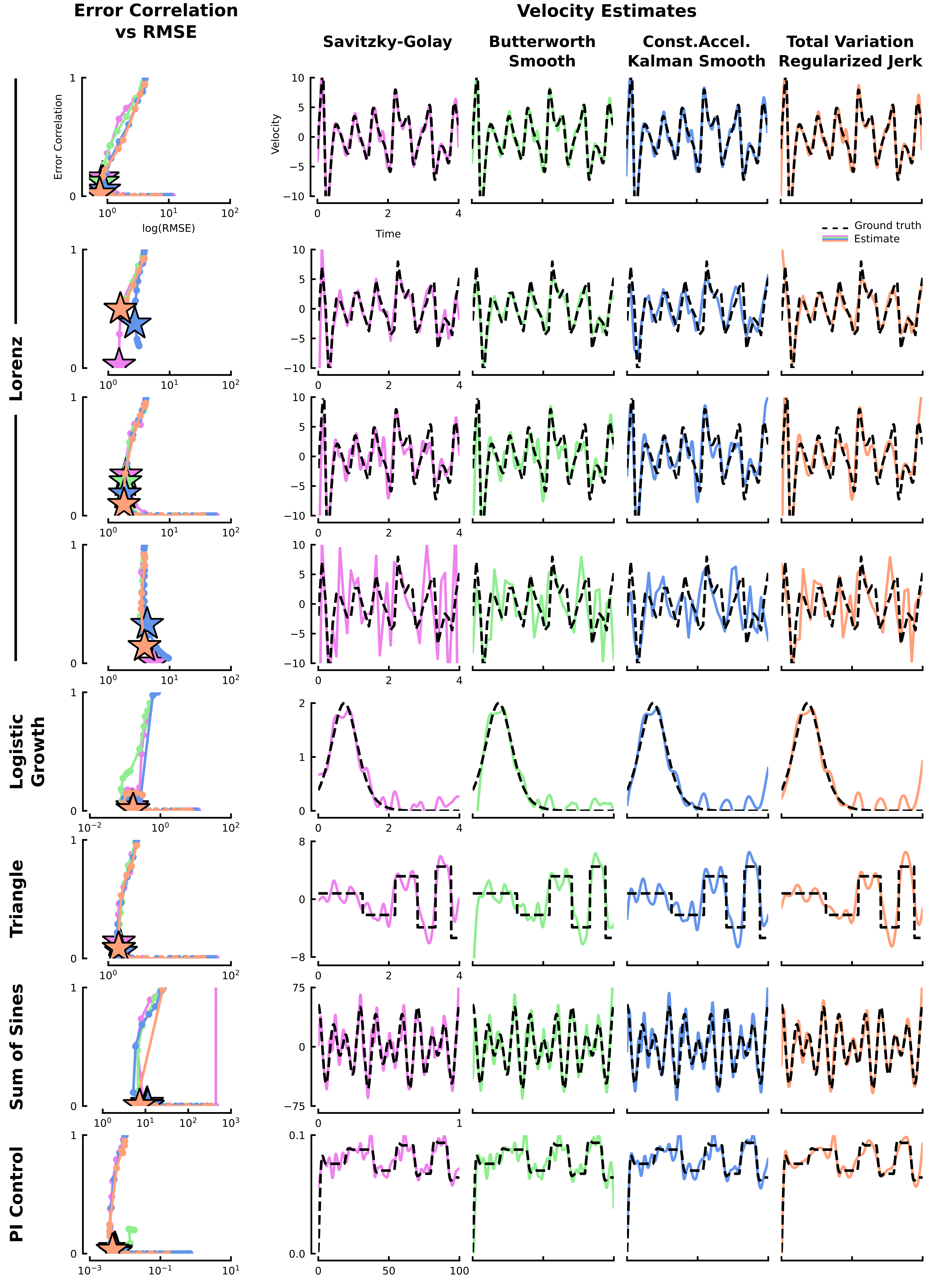}
\caption{Loss function and heuristic for choosing $\gamma$ is equally effective for different differentiation methods across a range of toy problems. Data plotted as in Fig. \ref{fig:fig_4}C, for each of the scenarios presented in Fig. \ref{fig:fig_3}. }
\label{fig:fig_5}
\end{figure*}

\subsection{Direct comparison of differentiation methods}

We examine how our loss function and heuristic for choosing $\gamma$ might perform on other differentiation methods beyond the Savitzky-Golay filter.
Figure~\ref{fig:fig_4} shows that for a noisy Lorenz system, the possible solution space is similar for all four methods we highlighted earlier, and our loss function achieves a similar Pareto front in each case. 
Note that although the Savitzky-Golay and Butterworth filters both operate in the frequency domain, the Kalman smoother and TVRJ methods do not. 

Interestingly, for all four differentiation methods, the possible solutions (the gray dots), and in particular their Pareto front, are quite similar, with the exception of the TVRJ method. 
This deviation may be because the TVRJ method only contains a single parameter. 
Our loss function, which defines the colored curves in the RMSE vs error correlation space, results in similar curves for each method, each of which follows the Pareto front quite closely. 
Although there are some differences in the location along the Pareto front that our heuristic selects as the optimal choice for each method, the resulting derivative estimates are qualitatively quite similar. 
A close comparison of the curves defined by the loss function, and the points selected by the heuristic, suggest that the Kalman and TVRJ methods produce slightly more accurate derivative estimates with a lower error correlation. 
However, looking at the resulting derivatives we see that the regions where the derivative estimates have high errors, all four estimates exhibit similar errors, suggesting that these errors may be a result of the data, not the method. 

These results suggest that our optimization framework is universal across different methods, a claim further supported by its performance across a range of synthetic problems (Fig.~\ref{fig:fig_5}.  
The most significant result of this analysis is that all four methods, despite being very different in their underlying mathematics, behave similarly under both our loss function and heuristic for choosing $\gamma$ across a wide range of data. 
Even in the case where they disagree on a quantitative level (second row, low temporal resolution Lorenz data), and the Savitzky-Golay filter appears to provide the estimate with the lowest error correlation, the resulting derivative estimates are in fact qualitatively quite similar. 

Taking a closer look at the errors in the derivative estimates across the range of toy problems shown in Fig.~\ref{fig:fig_5} reveals a subtle point about the limitations of the differentiation methods we highlight here. 
For all four methods, the errors in the derivative estimates are largest for the triangle problem, and to a lesser extent the proportional-integral control problem. 
These errors likely stem from two particular challenges. 
First, the frequency content of the data is very heterogeneous: it is near zero between the peaks and valleys, and near infinite at the peaks and valleys. Furthermore, the frequency of the oscillations for the triangle increase with time. 
Second, all four of the methods we highlighted here are designed to provide smooth derivatives, whereas the true derivative for the triangle problem is piece-wise constant. If this were known from the outset, it might be more effective to choose a method that is designed to return piece-wise constant derivatives, such as the total variation regularized on the $1^{st}$ derivative.

\section{Demonstrations on real-world data}

\begin{figure*}[htb] 
\includegraphics[width=0.9\textwidth]{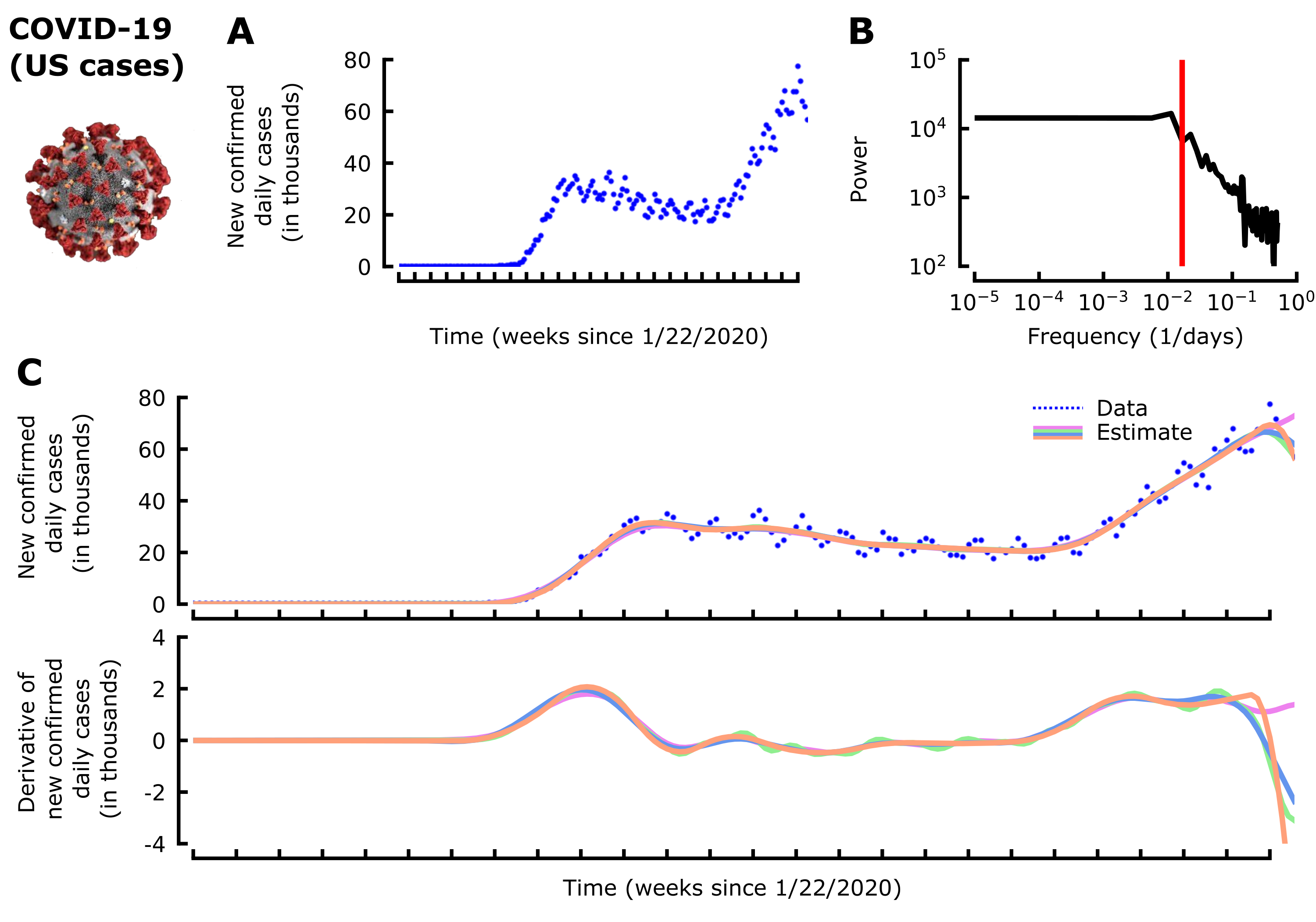}
\caption{Numerical differentiation of new confirmed daily cases in the United States of COVID-19 \cite{Dong2020} with no parameter tuning. A. Raw new daily cases. B. Power spectra of the data, indicating the cutoff frequency (red) used for selecting $\gamma=4.1$. C. Smoothed new daily cases, and their derivative, using a Savitzky-Golay filter (violet), a Butterworth filter (green), a constant acceleration Kalman forward-backward smoother (light blue), and total variation regularized jerk (orange). Note the similarity between all four methods except in the very last week, despite the significant differences in how each method works and the automated parameter selection.}
\label{fig:fig_6}
\end{figure*}

\begin{figure*}[p] 
\includegraphics[width=0.9\textwidth]{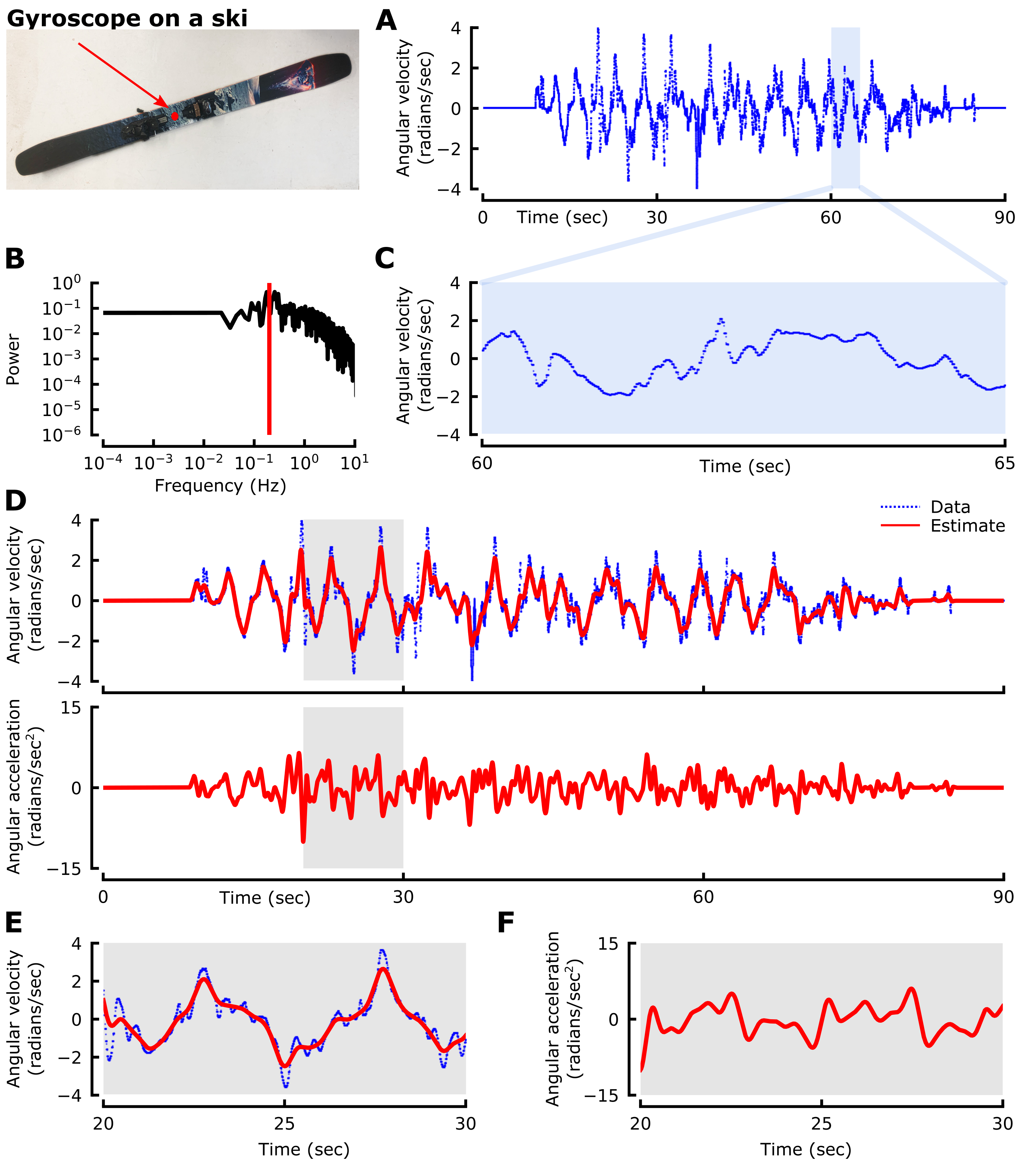}
\caption{Numerical differentiation of noisy gyroscope data from a downhill ski during one ski run, with no parameter tuning. A. Data from one axis of a gyroscope attached to the center of a downhill ski. B. Power spectra of the data, indicating the cutoff frequency (red) used for selecting $\gamma=11.5$. C. Zoomed in section of the data from A, which was used to optimize parameter selection. D. Smoothed angular velocities and angular accelerations, calculated using a Savitzky-Golay filter and the optimal parameters determined using our heuristic and loss function. E-F. Zoomed in sections from D.}
\label{fig:fig_7}
\end{figure*}

The real value of our multi-objective optimization framework is its straightforward application to real, noisy data where no  ground truth data is available. 
Here we provide two such examples: differentiation of the new confirmed daily cases in the United States of COVID-19, the disease caused by SARS-CoV-2 (Fig.~\ref{fig:fig_6}), and differentiation of gyroscope data from a downhill ski (Fig.~\ref{fig:fig_7}). 
In both examples, we examine the power spectra of the data to choose a cutoff frequency that corresponds to the start of the dropoff in power. 
This cutoff frequency, in conjunction with the time resolution of the data, are then used as inputs to our heuristic described by Eq.~\eqref{eqn:goldgamma} to determine an optimal value of $\gamma$. 
With $\gamma$ chosen, we minimize our loss function from Eq.~\eqref{eq:loss} to find the optimal parameters for numerical differentiation.

The year 2020 has seen a dramatic growth of the prevalence of a novel coronavirus, SARS-CoV-2, which causes the disease known as COVID-19. 
Estimating and understanding the rate of increase of disease incidence is important for guiding appropriate epidemiological, health, and economic policies.
In the raw data (~\cite{Dong2020}, https://github.com/CSSEGISandData/COVID-19) for the raw new confirmed daily cases of COVID-19, Fig. \ref{fig:fig_6}A) there is a clear oscillation with a period of one week, most likely due to interruptions in testing and reporting during weekends.
As such, we selected a lower cutoff frequency of 2 months, corresponding to the beginning of the steep drop off in the power spectra (Fig.~\ref{fig:fig_6}B). 
If the weekly oscillations were important, one could just as easily select a cutoff frequency of 1/week. 
Our heuristic for choosing $\gamma$ was based on sinusoidal data with a limited domain of time resolutions ranging from 0.001 to 0.1 seconds, so we scaled the time step units of the COVID-19 data to be close to this range, using $dt =$ 1 day, rather than 86,400 seconds. 
Our chosen cutoff frequency yielded a value of $\gamma=4.1$. 

Using this same value of $\gamma$ for each of the four differentiation methods under consideration resulted in very similar smoothed daily case estimates and derivatives, except during the final 2 weeks (Fig.~\ref{fig:fig_6}C). 
This highlights an important application of our method, which facilitates easy and fair comparison between different smoothing methods. 
Where these methods disagree, it is clear that none of the estimates can be trusted. 
A more subtle difference between the methods is that the Butterworth filter appears to preserve a larger remnant of the weekly oscillations seen in the raw data.

Finally, we consider angular velocity data collected from a gyroscope attached to a downhill ski over one minute of descent (Fig.~\ref{fig:fig_7}A) (ICM-20948, SparkFun; Wildcat Ski, Moment Skis). 
This type of data is representative of kinematic data that might be collected during experiments with robots or animals, which might be used to construct data-driven models of their dynamics~\cite{karashchuk2020anipose}. 
From the power spectrum, we chose a cutoff frequency of 0.2 Hz (Fig. \ref{fig:fig_7}B). 
This selection together with the time resolution of 0.0009 seconds yielded an optimal value of $\gamma=11.5$ using our heuristic. 
We calculated the smoothed angular velocity and acceleration estimates using a Savitzky-Golay filter (Fig. \ref{fig:fig_7}D--F). 
The other methods showed similar results (not shown for visual clarity), though the total variation method is not recommended for large datasets like this one due to the compounding computational costs.

\section{Discussion}

In summary, this paper develops a principled multi-objective optimization framework to provide clear guidance for solving the ill-posed problem of numerical differentiation of noisy data, with a particular focus on parameter selection. 
We define two independent metrics for quantifying the quality of a numerical derivative estimate of noisy data: the RMSE and error correlation. 
Unfortunately, neither metric can be evaluated without access to ground truth data. 
Instead, we show that the total variation of the derivative estimate, and the RMSE of its integral, serve as effective proxies.
We then introduced a novel loss function that balances these two proxies, reducing the number of parameters that must be chosen for any given numerical differentiation method to a single universal hyperparameter, which we call $\gamma$. 
Importantly, the derivative estimates resulting from a sweep of $\gamma$ lie close to the Pareto front of all possible solutions with respect to the true metrics of interest.
Although different applications may require different values of $\gamma$ to produce more smooth or less biased derivative estimates, we derive an empirical heuristic for determining a general purpose starting point for $\gamma$ given two features that can easily be determined from timeseries data: the cutoff frequency and time step. 
Our method also makes it possible to objectively compare the outputs for different methods. 
We found that for each problem that we tried, the four differentiation methods we explored in depth, including both local and global methods, 
all produce qualitatively similar results.

In our loss function we chose to use the RMSE of the integral of the derivative estimate and the total variation of the derivative estimate as our metrics. However, our loss function can be extended to a more general form,
\begin{equation}
     L = M_1 (\bxdothat, \bxdot) + \gamma_2 M_2 (\bxdothat, \bxdot) + \cdots + \gamma_p M_p (\bxdothat, \bxdot),
\end{equation}
where $M_1, M_2, \cdots , M_p$ represent $p$ different metrics that could be used, balanced by $p-1$ hyper-parameters. Alternative metrics include, for example, the tortuosity of the derivative estimate, the error correlation between the data and the integral of the derivative estimate, a metric describing the distribution of the error between the data and the integral of the derivative estimate. Depending on the qualities of the data and the specific application, different sets of metrics may be suitable as terms in the loss function.

Our loss function makes three important assumptions that future work may aim to relax. 
The first is that we assume the data has consistent zero-mean Gaussian measurement noise. 
How sensitive the loss function and heuristic are to outliers and other noise distributions remains an open question. 
It is possible that once we include other noise models, we will find differences in the behavior of differentiation methods. 
The second major limitation is that our loss function finds a single set of parameters for a given time series. 
For data where the frequency content dramatically shifts over time, it may be better to use time-varying parameters. 
Presently, this is limited by our current implementation, which relies on a computationally expensive optimization step. 
Future efforts may focus on ways to improve the efficiency of these calculations. 
Finally, we have focused on single dimensional time-series data. 
In principle, our proposed loss function can be used with multi-dimensional data, such as 2- and 3-dimensional spatial data, with only minor modifications. 

By simplifying the process of parameter selection for numerical differentiation to the selection of a single hyper-parameter, our approach makes it feasible to directly compare the performance of different methods within a given application.
One particular application of interest is that of data-driven model discovery. Methods such as sparse identification of nonlinear dynamics (SINDy)~\cite{Brunton2016pnas}, for example, rely directly on numerical derivative estimates, and the characteristics of these estimates can have an important impact on the resulting models. 
Using our method, it is now tractable to systematically investigate the collection of data-driven models learned from estimated derivatives of different smoothness and explore their impact on the models.

\section*{Acknowledgements}
We are grateful for helpful discussions with Steve Brunton and Pierre Karashchuk.
FvB acknowledges support from a Moore/Sloan Data Science and Washington Research Foundation Innovation in Data Science Postdoctoral Fellowship, a Sackler Scholarship in Biophysics, and the National Institute of General Medical Sciences of the National
Institutes of Health (P20GM103650). 
JNK acknowledges support from the Air Force Office of Scientific Research (FA9550-19-1-0011).
BWB acknowledges support from the Washington Research Foundation, the Air Force Office of Scientific Research (FA9550-19-1-0386), and the National Institute of Health (1R01MH117777).

\bibliography{references} 
\bibliographystyle{ieeetr}

\end{multicols}

\end{document}